\documentclass[12pt]{article}
\usepackage{amscd, amssymb, latexsym, amsmath}

\newtheorem{cor}{Corollary}
\newtheorem{pro}{Proposition}

\newenvironment{proof}
{\noindent {\em Proof}} {\hfill $\Box$}

\numberwithin{thm}{section}
\numberwithin{cor}{section}
\numberwithin{pro}{section}
\numberwithin{lem}{section}
\numberwithin{dfn}{section}
\numberwithin{rem}{section}
\numberwithin{equation}{section}

\newcommand{\R}{\mathbb R}
\newcommand{\C}{\mathbb C}

\begin{document}
\title{On Graphic Bernstein Type Results \\in Higher Codimension
}
\author{Mu-Tao Wang }
\date{January 31, 2002}
\maketitle
\vskip 10pt \centerline{email:\,
mtwang@math.columbia.edu}

\begin{abstract}
Let $\Sigma$ be a minimal submanifold of $\R^{n+m}$ that can be
represented as the graph of a smooth map $f:\R^n\mapsto\R^m$. We
apply a formula we derived in the study of mean curvature flow to
obtain conditions under which $\Sigma$ must be an affine subspace.
Our result covers all known ones in the general case. The
conditions are stated in terms of the singular values of $df$.
\end{abstract}

\section{Introduction}

The well-known Bernstein theorem asserts any complete minimal
surface that can be written as the graph of a function on $\R^2$
must be a plane. This result has been generalized to $\R^n$, for
$n\leq 7 $ and general dimension under various growth condition,
see \cite{eh} and the reference therein for the codimension one
case. In this note, we study the higher codimension case, i.e.
minimal submanifolds in $\R^{n+m}$ that can be written as graphs
of vector-valued functions $f:\R^n\mapsto \R^m$.

For higher codimension Bernstein type problems, there are general
results of \cite{fc}, \cite{hjw} and \cite{jx}. The idea in these
papers is to find a subharmonic function whose vanishing implies
$\Sigma$ is totally geodesic. Under the assumption that $f$ is of
bounded gradient, one can perform blow down analysis to reduce to
minimal cone case.  Maximum principle together with Allard
regularity theorem then complete the argument. A fundamental fact
is the Gauss map of a minimal submanifold is a harmonic map, so
any convex function on the Grassmannian renders a subharmonic
function. The work in \cite{fc}, \cite{hjw}, and \cite{jx}
consists of delicate analysis of the geometry of Grassmannian in
order to locate the maximal region where a convex function exists
. The condition for Bernstein type reuslts is in terms of the
function

\[*\Omega_1=\frac{1}{\sqrt{\det(I+(df)^t df)}}\]
This can be considered as the inner product of  the tangent space
of $\Sigma$ and the domain $\R^n$. In particular, in \cite{fc} and
\cite{hjw}, it was shown,

\[*\Omega_1\geq K'> \cos^p(\frac{\pi}{2\sqrt{2}p})\]
where $p=\min(n,m)$ implies $\Sigma$ is an affine subspace. This
bound is later improved in \cite{jx} to

\[*\Omega_1\geq K''> \frac{1}{2}\]

We adopt a direct approach to calculate explicitly the Laplacian
of $\ln *\Omega$. The calculation first appeared to us in
\cite{mu3} in the context of mean curvature flow. The formula
determines precisely when $\ln *\Omega$ is a superharmonic
function.

\vskip 10pt \noindent {\bf Theorem A} {\it Let $\Sigma$ be a
minimal submanifold of $\R^{n+m}$. Suppose $\Sigma$ can be written
as the graph of a smooth map $f:\R^n\mapsto\R^m$. Denote the
singular values of $df$ by $\lambda_i, i=1,\cdots, n$. If there
exist $\delta, K>0$ such that $|\lambda_i\lambda_j|\leq 1-\delta$
for any $i\not= j$ and $*\Omega_1=\frac{1}{\sqrt{\det(I+(df)^t
df)}}\geq K$, then $\Sigma$ is an affine subspace. } \vskip 10pt

Notice that $K$ is not necessarily bounded from below. At the end
of \S3, we discuss how Theorem A implies the results in \cite{fc},
\cite{hjw} or \cite{jx}.

We remark Theorem A also holds for submanifolds with parallel mean
curvature vector. This condition for minimal submanifolds is not
the optimal one. We state it first because it is simple to apply
and seems natural ( see the remark at the end of the introduction
on the geometric meaning). The optimal condition, which is in
terms of $\lambda_i$ again, is somewhat more complicate. It will
be described explicitly in \S 2.

\vskip 10pt \noindent {\bf Theorem B} {\it Under the assumption of
Theorem A. There exists an optimal condition in terms of the
singular values of $df$ such that $\Sigma$ is an affine subspace
whenever $f$ satisfies this optimal condition } \vskip 10pt

Theorem B also implies the result of \cite{b} and \cite{fc} that
any three-dimensional minimal cone is flat and thus is sharper
than Theorem A.

The higher codimension Bernstein type result is not expected to be
true in the most generality due to the counterexample of
Lawson-Ossermann \cite{lo}.

Notice that in codimension one case, any graphic minimal
hypersurface is stable. It is thus natural to impose stable
assumption in higher codimension. A special category of minimal
submanifolds are calibrated submanifolds studied by Harvey and
Lawson \cite{hl}. They are volume minimizing and thus stable.
These include holomorphic submanifolds, special Lagrangian
submanifolds, assoicative, coassociative, and Cayley submanifold.
The special Lagrangian case  is studied in \cite{fu},\cite{n}
\cite{jx2} \cite{y}, and \cite{tw}. The Lagrangian condition
implies the existence of a potential function $F$ such that
$f=\nabla F$. In terms of PDE, the Bernstein type problem asks
when an entire solution of

\begin{equation}\label{sp}
 Im(\det(I+\sqrt{-1} D^2F))=constant
 \end{equation}
is a quadratic polynomial. In this case, the condition of Theorem
A or B can be stated in terms of the eigenvalues of $D^2F$.

We remark that conditions on the singular values of $df$ or
eigenvalues of $D^2F$ depend on the choice of a subspace (usually
a calibrated one) over which $\Sigma$ is written as a graph. Such
conditions correspond to regions in the relevant Grassmannian. A
more invariant description is to require the image of the Gauss
map lies in the orbit of one of these regions under the respective
holonomy group actions.

To demonstrate, let us compare the condition
$|\lambda_1\lambda_2|<1$ and $(1+\lambda_1^2)(1+\lambda_2^2)<4$
\,(or
$*\Omega_1=\frac{1}{\sqrt{(1+\lambda_1^2)(1+\lambda_2^2)}}>\frac{1}{2}$)
when $n=m=2$. It is clear the former condition is weaker than the
latter one. In this case, $G(2,4)=S^2(\frac{1}{\sqrt{2}})\times
S^2(\frac{1}{\sqrt{2}})$. In $f:\R^2 \mapsto \R^2$, denote the
coordinates on the domain $\R^2 $ by $x^1, x^2$ and the
coordinates on the target $\R^2 $ by $y^1, y^2$. Then the forms
$\omega_1, \omega_2$

\[\omega_1=\frac{1}{\sqrt{2}}(dx^1\wedge dx^2-dy^1\wedge dy^2)\]

\[\omega_1=\frac{1}{\sqrt{2}}(dx^1\wedge dx^2+dy^1\wedge dy^2)\]
serve as coordinate functions on $G(2,4)$, see for example section
3 in \cite{mu1}. We may consider them as the height function on
each of the two $S^2(\frac{1}{\sqrt{2}})$. The condition
$|\lambda_1\lambda_2|<1$ corresponds to $*\omega_1>0$ and
$*\omega_2>0$, which means the image of the Gauss map lies in the
product of the two hemispheres. These are natural geometric
condition. It is shown in \cite{mu1} that under this assumption,
the maximum principle not only works in the elliptic case but also
in the parabolic one.

 Contrary to the codimension one case,
the equation of the second fundamental form has not been used in
higher codimension to our knowledge . We expect a complete
classification theorem of graphic calibrated submanifolds will
rely on such estimate.

The author wish to thank Professor D. H. Phong and Professor S.-T.
Yau for their  encouragement and  support. He also thanks Mao-Pei
Tsui for useful discussions.

\section{ General results}
We first recall a formula whose parabolic version was derived in
\cite{mu3}. To apply the formula in \cite{mu3} to the current
situation, we note that a minimal submanifold corresponds to
stationary phase of mean curvature flow.

Let $\Sigma$ be an $n$-dimensional submanifold of $\R^{n+m}$ and
$\Omega$ a parallel $n$ form on $\R^{n+m}$.
 Around any point $p\in \Sigma$, we choose any orthonormal frames
$\{e_i\}_{i=1\cdots n}$ for $T_p\Sigma$ and
$\{e_\alpha\}_{\alpha=n+1, \cdots, n+m}$ for $N_p\Sigma$, the
normal bundle of $\Sigma$ . The convention that $i, j, k, \cdots$
denote tangent indexes and $\alpha, \beta, \gamma \cdots $ denote
normal indexes is followed.

The second fundamental form of $\Sigma$ is denoted by  $h_{\alpha
ij} =<\nabla_{e_i} e_j, e_\alpha>$.

If we assume $\Sigma$ has parallel mean curvature vector, then
$*\Omega =\Omega(e_1, \cdots, e_n)$ satisfies (see Proposition 3.1
in \cite{mu3})

\begin{equation}\label{main}
\begin{split}
&\Delta*\Omega +*\Omega (\sum_{\alpha ,l,k}h_{\alpha l
k}^2)\\
&-2\sum_{\alpha, \beta, k}[\Omega_{\alpha \beta 3\cdots n}
h_{\alpha 1k}h_{\beta 2k} + \Omega_{\alpha 2
\beta\cdots,n}h_{\alpha 1k}h_{\beta 3k} +\cdots +\Omega_{1\cdots
(n-2) \alpha \beta} h_{\alpha (n-1)k}
h_{\beta nk}]=0\\
\end{split}
\end{equation}

\begin{equation}\label{gradient}
(*\Omega)_k =\sum_\alpha \Omega_{\alpha 2 \cdots n}h_{\alpha 1 k}
+\cdots +\Omega_{1 \cdots {n-1} \alpha} h_{\alpha n k}
\end{equation}
where $\Delta$ is the Laplace operator of the induced metric on
$\Sigma$ and $\Omega_{\alpha \beta 3\cdots n}=\Omega(e_\alpha,
e_\beta, e_3, \cdots, e_n)$.

Now take $\Omega=dx^1\wedge\cdots\wedge dx^n$, i.e the volume form
of $\R^n$. In this case $*\Omega=\Omega(e_1, \cdots, e_n)$ is
exactly the Jacobian of the projection from $\Sigma$ to $\R^n$. A
similar formula in this case appeared in Fischer-Colbrie's paper
\cite{fc}\,(Lemma 1.1).

When $m=1$, i.e. the codimension one case,
$\partial_i=\frac{\partial}{\partial x^i}+\frac{\partial
f}{\partial x^i} \frac{\partial}{\partial x^{n+1}}$ forms a basis
for $T_p\Sigma$. Then

\[*\Omega=\frac{\Omega(\partial_1,\cdots,
\partial_n)}{\sqrt{\det(\partial_i,\partial_j)}}\]

It is not hard to compute that
$*\Omega=\frac{1}{\sqrt{1+|df|^2}}$, i.e. the angle made by the
normal vector of $\Sigma$ and the $x_{n+1}$ axis. Thus formula
(\ref{main}) reduces to the classical formula for nonparametric
minimal surfaces.

\[\Delta
\frac{1}{\sqrt{1+|df|^2}}+\frac{1}{\sqrt{1+|df|^2}}|A|^2=0\]

A standard argument will imply the Bernstein type result for
bounded $|df|$. In fact, this equation coupled with Simon's
identity gives the optimal results in codimension one case, see
\cite{eh}.

In the general dimensional case, we choose a particular basis to
represent formula (\ref{main}) for $\Omega=dx^1\wedge\cdots\wedge
dx^n$. Given any $p$ then the differential of $f$ is a linear map
from $\R^n$ to $\R^m$. We can use singular value decomposition to
find orthonormal bases $\{a_i\}_{i=1\cdots n}$ for $\R^n$ and
$\{a_\alpha\}_{\alpha=n+1,\cdots, n+m}$ for $\R^m$ such that

\[df(a_i)=\lambda_i a_{n+i}\]
for $i=1\cdots n$. Notice that $\lambda_i=0$ if $i>\min\{n, m\}$.

 Now $\{e_i= \frac{1}{\sqrt{1+\lambda_{i}^2}}
(a_i+\lambda_{i}a_{n+i})\}_{i=1, \cdots, n}$ forms an orthonormal
basis for $T_p\Sigma$ and $\{e_{\alpha}=
 \frac{1}{\sqrt{1+ \lambda_{\alpha-n}^2}}
(a_{\alpha}-\lambda_{\alpha -n}a_{\alpha-n})\}_{\alpha=n+1,\cdots,
n+m} $ an orthonormal basis for $N_p\Sigma$.

It is not hard to see $*\Omega=\frac{1}
{\sqrt{\prod_{i=1}^n(1+\lambda_{i}^2)}}$. Apply these bases to
equation (\ref{main}).

\begin{pro}\label{equation}
Let $\Sigma$ be the graph of a smooth map $f:\R^n\mapsto \R^m$ and
$\{\lambda_i\}_{i=1\cdots n}$ be the singular values of $df$. If
the mean curvature vector of $\Sigma$ is parallel,  then $*\Omega$
satisfies the following equation.
\begin{equation}\label{eq}
\begin{split}
&\Delta *\Omega= -*\Omega\{\sum_{\alpha, l, k} h_{\alpha
lk}^2-2\sum_{k, i<j} \lambda_{ i}\lambda_{ j}h_{n+i,ik} h_{n+j,
jk} +2\sum_{k, i<j}\lambda_{i}\lambda_{j} h_{n+j, ik} h_{n+i,
jk}\}
\end{split}
\end{equation}
where $\Delta $ is the Laplace operator of the induced metric on
$\Sigma$. The indexes $i,j$ range between $1$ and $\min\{n, m\}$.

\end{pro}

\vskip 10pt
\begin{proof}\textit{\,of Theorem A.}

We shall calculate
\begin{equation}\label{eq2}
\Delta ( \ln *\Omega ) = \frac{ *\Omega \Delta ( *\Omega ) -
|\nabla *\Omega|^2}{ |*\Omega|^2 }
\end{equation}

By formula (\ref{gradient}), the covariant derivative of $*\Omega$
is .

 \[(*\Omega)_{k} = -*\Omega
(\sum_{i} \lambda_{i}h_{n+i,ik} ) \]

Plug this and equation (\ref{equation}) into equation (\ref{eq2})
and we obtain

\begin{equation}\label{same}
\Delta ( \ln *\Omega )= - \{\sum_{ \alpha,l, k} h_{\alpha lk}^2 +
\sum_{k, i} \lambda_{ i}^2 h_{n+i, ik}^2 +2\sum_{k,
i<j}\lambda_{i}\lambda_{j} h_{ n+i,jk}h_{n+j, ik} \}
\end{equation}

From the assumption of the theorem, it is obvious that $\Delta (
\ln *\Omega )\leq -\delta_1|A|^2$ by completing square. The
condition $\frac{1}{\sqrt{\det(I+(df)^t df)}}\geq K$ means
$\Sigma$ is the graph of a vector-valued function with bounded
gradient. Therefore we can perform blow down to get a minimal
cone. Notice that these conditions are invariant under scaling. We
can apply maximum principle to conclude the minimal cone is flat
and then Allard regularity theorem forces $\Sigma$ to be an affine
space.

\end{proof}

\vskip 10pt
\begin{proof}\textit{\,of Theorem B.}
The optimal condition can be found by considering

\[F(h_{\alpha lk})=\sum_{ \alpha,l, k} h_{\alpha lk}^2 +
\sum_{k, i} \lambda_{ i}^2 h_{n+i, ik}^2 +2\sum_{k,
i<j}\lambda_{i}\lambda_{j} h_{ n+i,jk}h_{n+j, ik}\] as a quadratic
form define on the space of all possible $h_{\alpha lk}$ with
$h_{\alpha lk}=h_{\alpha k l}$ and $\sum_k h_{\alpha kk}=0$ for
each $\alpha$. Now the problem is reduced to determine conditions
on $\lambda_i\lambda_j$ to guarantee $F$ is a positive quadratic
form with $F(h_{\alpha lk})\geq \epsilon \sum h^2_{\alpha lk}$ for
some $\epsilon>0$. It is computable, yet it seems hard to find a
simple description. At the end of next section, we show how this
gives the optimal result when $n=2$.

\end{proof}

The condition in Theorem A or B depends on the choice a subspace
over which $\Sigma$ is written a graph. Since being minimal is
invariant under the isometry action of $O(n+m)$, a slightly more
general condition is the following.

\begin{cor}
If $\Sigma$ is the graph of $f:\R^n\mapsto \R^m$. Denote $A=df$ as
an $n\times m$ matrix . If there exists an element $ g\in O(n+m)$
of the block form $g= \left[
\begin{matrix} P&Q
 \\ R&S \end{matrix} \right]$  such that $P+AR$ is invertible and
 the singular values of $(P+AR)^{-1}(Q+AS)$ satisfy the optimal condition in Theorem B
   , then $\Sigma$
 is an affine subspace.
\end{cor}
\begin{proof}
 The tangent space of $\Sigma$ is represented by
$\left[\begin{matrix}I&A\end{matrix}\right]$. Under the action of
$g$ on $\R^{n+m}$, the tangent space of $g(\Sigma)$ becomes
$\left[\begin{matrix}P+AR&Q+AS\end{matrix}\right]$. If $P+AR$ is
invertible, then $g(\Sigma)$ is still the graph of some $\bar{f}$
over $\R^n$. Now $d\bar{f}=(P+AR)^{-1}(Q+AS)$.
\end{proof}

\vskip 10pt
 When $\Sigma$ is a Lagrangian submanifold defined as
the graph of $\nabla F:\R^n\mapsto \R^n$.  Let $A= D^2 F $. Note
being minimal Lagrangian is invariant under $U(n)$. An element $g$
in $U(n)$ can be expressed as a $2n \times 2n$ block matrix of the
form

\[\left[ \begin{matrix} P&-Q
 \\ Q&P \end{matrix} \right]\]
 with
 \[PP^t+QQ^t=I, -PQ^t+QP^t=0\]
where the first $n$ components correspond to the real $x^i$ and
the last $n$ components the imaginary $y^i$.

Rotate $\C^n$ by this element $g$, the tangent space  becomes
$\left[
\begin{matrix} P+AQ & -Q+AP \end{matrix} \right]$
If we require $P+AQ$ is invertible then the condition is
determined by the eigenvalues of the symmetric matrix
$(P+AQ)^{-1}(-Q+AP)$. Please see \cite{tw} for the optimal
condition in this case.

\section{Applications}

In the last section, we show how our theorems imply all known
general graphic Bernstein type results in higher codimension. In
Theorem 1 of \cite{jx}, the authors proved a Bernstein type
theorem in terms of a bound on

\[\Delta_f=\{\det(\delta_{\alpha \beta}+\sum_i f^i_{x_\alpha}(x)
f^i_{x_\beta}(x))\}^{\frac{1}{2}}<2\] This result improves those
in \cite{fc} and \cite{hjw}. This quantities of course is our
$\frac{1}{*\Omega}={\sqrt{\prod_{i=1}^n(1+\lambda_{i}^2)}}$.

It is not hard to check that $\prod_{i=1}^n(1+\lambda_{i}^2)<4$
implies $|\lambda_i\lambda_j|<1$ for $i\not=j$. Therefore, Theorem
A implies the result in \cite{jx}.

As was mentioned above, Theorem A only assumes $\Sigma$ is of
parallel mean curvature vector. The condition is weakened if we
further assume $\Sigma$ is minimal, i.e $\sum_k h_{\alpha kk}=0$.
We demonstrate how it works in two-dimension. Write out the right
hand side of equation (\ref{same}) when $n=2$.

\begin{equation}
|A|^2 + \lambda_1^2(h_{n+1,1,1}^2+h_{n+1,
1,2}^2)+\lambda_2^2(h_{n+2,2,1}^2+h_{n+2,
2,2}^2)+2\lambda_1\lambda_2(h_{n+1,
2,1}h_{n+2,1,1}+h_{n+1,2,2}h_{n+2,1,2})
\end{equation}

Since $\Sigma$ is minimal, $h_{n+1,1,1}=-h_{n+1, 2,2}$ and
$h_{n+2,2,2}=-h_{n+2,1,1}$. This can be completed square and
becomes

\[|A|^2+(\lambda_1 h_{n+1, 2,2}+\lambda_2
h_{n+2,1,2})^2+(\lambda_1 h_{n+1, 1,2}+\lambda_2 h_{n+2,1,1})^2\]

This in fact applies to three dimensional minimal cone because the
second fundamental form vanishes in one direction and the right
hand side of equation (\ref{same}) reduces to the two-dimensional
case. This recovers the result of Barbosa \cite{b} and
Fischer-Colbrie \cite{fc}.

\end{document}